\magnification=\magstep1
\parindent 0cm

\def\Bbb R{I\! \! R}
\def\R{I\! \! R}
\def \bull{\vrule height .9ex width .8ex depth -.1ex}
\def\P {I\! \! P}
\def\pj {P_j(x,\xi)}
\def\pjs {P_j^{2}(x,\xi^*)}
\def\kj {k_j(\xi)}
\def\gj {\gamma_j(\xi^*)}
\def\kji {k_j^{-1}(\xi^*)}

\def\kn {k_n(\xi)}

\pageno=1

\smallskip
\centerline{A NOTE ON SOME PECULIAR NONLINEAR EXTREMAL}
\centerline{PHENOMENA OF THE CHEBYSHEV POLYNOMIALS}
\centerline{by}
\smallskip
\centerline{Holger Dette$^*$}  
\centerline{Institut f\"ur Mathematische Stochastik}  
\centerline{Abteilung Mathematik}
\centerline{Technische Universit\"at Dresden}
\centerline{Mommsenstr.\ 13}
\centerline{01062 Dresden}
\centerline{GERMANY}
\footnote{}{AMS Subject Classification: 33C45\hfill\break
		  \hfill\break
		  $^*$Research supported in part by the Deutsche
		  Forschungsgemeinschaft\hfill\break
		  }
\bigskip
\bigskip
\centerline{ABSTRACT}

\medskip
{
We consider the problem of maximizing the sum of squares of the leading
coefficients of polynomials $P_{i_1}(x),\ldots ,P_{i_m}(x)$ (where $P_j(x)$
is a polynomial of degree $j$) under the restriction that the sup-norm
of $\sum_{j=1}^m P_{i_j}^2(x)$ is bounded on the interval $[-b,b]$ ($b>0$).
A complete solution of the problem is presented using duality theory of
convex analysis and the theory of canonical moments. It turns out,
that contrary to many other extremal problems the structure of the solution
will depend heavily on the size of the interval $[-b,b]$.}

\bigskip
\bigskip

\baselineskip=14.6pt

\noindent {\bf 1.\ \ Introduction.}
Let $\P_j$ denote the set of all polynomials of degree $j$, $I=
\{i_1,\ldots ,i_m\}$
denote a subset of $\{1,\ldots ,n\}$ containing $n$ (i.e. $n\in I$,
$i_l \neq i_k$ if $k \neq l$) and define
$$
P_I~:=~\Bigl\{ (P_j)_{j\in I}~|~P_j \in \P_j, j\in I, ~
\sup_{x\in [-b,b]} \sum_{j\in I} P_j^2(x) \leq 1\Bigr\}
$$
as the set of all polynomials of degree $i_1,\ldots ,i_m$ such that
the sup-norm of the sum of squares is bounded by $1$ 
on the interval $[-b,b]$.
 In the following $m_l(P_l)$ denotes the leading coefficient of
the polynomial $P_l \in \P_l$ and we are interested in the nonlinear
extremal problem
$$
\max \Bigl\{ \sum_{l\in I} m_l^2(P_l) ~|~(P_l)_{l\in I} \in P_I \Bigr\}~~.
\leqno ({\cal P}_I)
$$
For $I=\{n\}$ $({\cal P}_I)$ yields the well known extremal property of the
Chebyshev polynomials of the first kind $T_n({x\over b})$ (see e.g. 
Natanson (1955), Achieser (1956) or Rivlin (1990)). Similar problems 
were investigated by Dette (1994a) who considered the maximization of a weighted 
product of the squared leading coefficients of the polynomials 
$P_l(x)$.  All extremal problems in these references satisfy
 a so called
``invariance property'' which means that if a solution on one interval, 
say $[-1,1]$, has been determined, then the solution  on another interval can 
easily be obtained by a linear transformation from the ``optimal'' polynomials
on the interval $[-1,1]$.

In this note we will present a complete solution of the (nonlinear)
extremal problem $({\cal P}_I)$. It will turn out that the above 
invariance property is not true any longer for the problem $({\cal P}_I)$
if $m \ge 2$.  While for sufficently small $b>0$ the Chebyshev polynomial
(of the first kind)
on the interval $[-b,b]$ of degree $\max_{j=1}^m i_j$ is a solution of
(${\cal P}_I$) (all other polynomials are vanishing) this is not true
any longer for large $b$.
Here the structure of the extremal solution depends heavily on the size of
the interval $[-b,b]$. 

In Section 2 the problem $({\cal P}_I)$ is solved
by an application of some results in convex analysis (see Pukelsheim (1993))
and the theory of canonical moments (see Studden (1981)). 
It turns out that the problem (${\cal P}_I$) is dual to a maximization
problem of a concave function defined on the set of all probability
measures on the interval $[-b,b]$. This problem appears in
the theory of optimal experimental design in
mathematical statistics (see Dette (1994)).
 While from a statistical point of view the support 
points and weigths of the optimal measure 
are the main interest it is shown in this paper
that the orthogonal polynomials with respect to
this measure form essentially the solution of the extremal problem
(${\cal P}_I$). Section 3
deals with some special cases for the set $I$, namely $I=\{1,\ldots ,n\}$
and $I=\{n-1,n\}$ and some explicit examples. Finally, in Section 4, 
similar problems are investigated which generalize the extremal properties
of the Chebyshev polynomials of the second kind.

\bigskip
\bigskip

{\bf 2. The Solution of (${\cal P}_I$).}~~Throughout this paper $\xi$
is a probability measure on the interval $[-b,b]$ and the corresponding 
orthogonal polynomials with leading coefficient $1$ will be denoted
by $\pj$ and their (squared) 
 $L_2$-norm by $\kj =\int_{-b}^b P_j^2(x,\xi) d\xi (x)$.
The main step for solving the extremal problem $({\cal P}_I)$
is the following duality which is proved in the appendix.

\bigskip

{\bf Theorem 2.1.}~~{\it Let $\Xi :=\{ \xi~|~ \kn >0\}$ and $n\in I$, then
$$
({\cal P}_I)~~~~~~~~~~~~
\max\Bigl\{ \sum_{l\in I}m_l^2(P_l)~|~(P_j)_{j\in I} ~\in ~P_I \Bigr\}
~=~\min_\xi \max \left\{ k_j^{-1}(\xi)~|~j\in I\right\}~~~~~~~~~~~
\eqno ({\cal D}_I)
$$
Moreover, if $\xi^*$ is a solution of the problem $({\cal D}_I)$
and 
$$
{\cal M}(\xi^*) = \{ j\in I~|~ k_j(\xi^*) = \min_{i\in I} k_i
(\xi^*) \},
$$
then $\{\sqrt{\alpha_j  /k_j(\xi^*)} P_j(x,\xi^*) \}_{j\in I}$ is a solution of $
({\cal P}_I)$ where $P_j(x,\xi^*)$ is the $j$th monic orthogonal
polynomial with respect to the measure $d\xi^* (x)$ and the $\alpha_j$ 
are (arbitrary) nonnegative numbers with sum $1$ satisfying
$$
\alpha_j=0 ~~~~~\hbox{if}~~j\in I\setminus {\cal M}(\xi^*)
\leqno (2.1)
$$
$$
\sum_{j\in I} \alpha_j \kji \pjs ~\leq ~1 ~~~\hbox{for all}
~~x\in [-b,b]~.
\leqno (2.2)
$$
}

\bigskip

The dual problem (${\cal D}_I$) appears in the theory of
optimal experimental
design in mathe- matical statistics and has been
solved in the special case $I=\{1,\ldots ,n\}$ (see Dette (1994)). 
While  statisticians are mainly interested into the support points and
weights of the solution $\xi^*$ of  (${\cal D}_I$) 
(they give essentially the points where observations have to
be taken in a polynomial regression)  Theorem 2.1 shows 
that the orthogonal polynomials with respect to the measure $d\xi ^*
(x)$ are needed for the solution of the primal problem
 (${\cal P}_I$). 
In order to  determine these polynomials (and 
to solve the dual problem $({\cal D}_I)$) some basic
facts about canonical moments of probability measures on the interval 
$[-b,b]$ are needed. The Stieltjes transform with 
corresponding continued fraction
 expansions of such a measure $\xi$ is given by
$$
\eqalign{
\int_{-b}^b {d\xi (x) \over z-x}~ =&~
{~ ~1~~~|\over |z+b}
~-~ {~2b\zeta_1~|\over \! \! \! \! \! \! \!|~ ~1} 
~-~ {~2b\zeta_2~|\over  |\ z+b} 
~-~ {~2b\zeta_3~|\over \! \! \! \! \! \! \! |~ ~1} 
~-~ {~2b\zeta_4~|\over  |\ z+b} 
~-~~\ldots \cr
=&~{\ \ \ 1\ \ \ \ 
~~~~~~~~~|\over |z+b
(1-2\zeta _1)}
~-~ {~~~~~(2b)^2\zeta_1\zeta_2~~~~~~~~~|\over  
|\ z+b(1-2\zeta_2-2\zeta_3)} ~-~ {~~~~~(2b)^2\zeta_3\zeta_4~~~~~~~~~|\over  
|\ z+b(1-2\zeta_4-2\zeta_5)} ~-~~\ldots \cr
}
$$
where $\zeta_1=p_1$, $\zeta_j=q_{j-1}p_j$ ($j\ge 2$), $q_j=1-p_j$
($j \geq 1$)
and $0\leq p_j \leq 1$ (see e.g. Lau and Studden (1988)).
The quantities $p_j$ are called the canonical moments of
$\xi$. Note that $p_{j+1}$ is undefined whenever $p_{j}
\in \{0,1\}$ because in this case the continued fraction terminates. It is well
known that the polynomial in the denominator of the $n$th convergent is the
$n$th  monic orthogonal polynomial with respect to the measure
$d\xi (x)$ and that this polynomial is given by
the continuant (see Perron, Bd. I, (1954), p. 9)
$$
\eqalign{
& (2.4) ~~~~~~~~ P_n (x,\xi)~=~ \cr
& \cr 
& K  \pmatrix{~~-(2b)^2 \zeta_1\zeta_2~~~~~~~~~~~~~~
-(2b)^2 \zeta_3\zeta_4 & \ldots& 
\! \! \! \! \! \! \! \! \! \! \! \! \! \! \! \!   -(2b)^2 \zeta_{2n-3}\zeta_{2n-2}\cr
            x+b(1-2\zeta_1)~~~x+b(1-2\zeta_2-2\zeta_3)
            \qquad ~~ & \ldots \ldots & ~~~x +b(1-2\zeta_{2n-2}
            -2\zeta_{2n-1})\cr} \cr
}
$$
and has $L_2$-norm
$$
k_n(\xi)~=~\int_{-b}^b P_n^2(x,\xi) d\xi (x) ~=~
(2b)^{2n} \prod_{j=1}^n \zeta_{2j-1}\zeta_{2j}
\leqno (2.5)
$$
(see Chihara (1978) or Wall (1948)). The following
theorem determines the canonical moments of the solution of
the dual problem $({\cal D}_I)$.

\bigskip

{\bf Theorem 2.2.}{\it ~~The solution $\xi^*$
of the dual problem $({\cal D}_I)$ 
is uniquely determined by its canonical moments $p_{2j-1}^*={1\over 2}$
$(j=1,\ldots ,n)$, $p_{2n}^*=1$ and 
$$
p_{2 (n-j)}^* = \max \Bigl\{ z_{n-j}\Bigl[1 - 
 b^{-2j} \prod_{i=
   n-j+1}^{n-1} ( q_{2i}^* p_{2i}^*
)^{-1}\Bigr], {1\over 2}\Bigr\}
~~~~~~~~~~~j=n-1,\ldots ,1. \leqno (2.6)
$$
where $z_{n-j}$ is $1$ or $0$ according to $n-j \in I$ or $n-j 
\notin I$.
}

\bigskip

{\bf Proof.}~~In the special case $I=\{1,\ldots ,n\}$ a proof of Theorem 2.2
can be found in Dette (1994), which can be generalized to
arbitrary index sets. For the sake of completeness
we provide a different proof in this paper, which
is directly based on the duality result of Theorem 2.1
and uses some identities for orthogonal polynomials
on campact intervals. Let $\gamma_{n-j} =
1 -  b^{-2j} \prod_{i=
   n-j+1}^{n-1} \left( q_{2i}^* p_{2i}^*\right)^{-1} $
($\gamma_n=1$), then
   it is easy to see (observing (2.5) and (2.6))
   that $\gamma_{n-j} \ge {1\over 2}
   $ if and only if $k_n(\xi ^*)
   = k_{n-j}(\xi ^*)$ and $\gamma_{n-j} <{1\over 2}
   $ if and only if $k_n(\xi ^*)< k_{n-j}(\xi ^*)$ ($n-j \in I$).
 Consequently
   we have for the set
   ${\cal M }(\xi ^*)$ in Theorem 2.1 and the canonical moments
defined in (2.6)
     $$
\leqalignno{
n \in {\cal M}(\xi^*) &  ~=~\{j\in I~|~\gamma_j \ge {1\over 2}\}
& (2.7) \cr
p_{2j}^*&~={1\over 2} ~~~\hbox{if } j\notin {\cal M}(\xi^*) .& (2.8)\cr
}
$$
In the following let $P_l(x,\xi^*)$ denote the
$l$th monic orthogonal polynomial with respect to the
measure $d\xi^*(x)$ and  define
$$
\alpha_j ~=~ \prod_{i=1}^{j-1} {q_{2i}^*\over p_{2i}^*} \left(1-
 {q_{2j}^*\over p_{2j}^*}\right),
 \leqno (2.9)
 $$
which have sum $1$ and are nonnegative, by the definition of $p_{2j}^*$ in
(2.6). From
 Theorem 3.5 and 4.1 in Dette (1994b) it follows that
the orthonormal
 polynomials $k_l^{-1/2}(\xi ^*)P_l(x,\xi ^*) $
 with respect to the measure $d\xi^*(x)$ satisfy
 $$
 \sum_{j=1}^n \alpha_j k_j^{-1}(\xi ^* )P_j^2(x,\xi ^*)~=~
 \sum_{j\in {\cal M}(\xi ^*)}\alpha_j k_j^{-1}(\xi ^*) P_j^2(x,\xi ^*)
~\leq ~1
\leqno (2.10)
$$
for all $x\in [-b,b]$. Note that the result in Dette (1994b)
was originally stated on the interval $[-1,1]$ but can easily
be transferred to the interval $[-b,b]$ and that we
have  used $\alpha _j \geq 0$,
$\alpha_j =0$ if  $j\notin {\cal M}(\xi^*)$, which follows from (2.8) 
and (2.9). By (2.10) we thus have
$$
\{P_j^*(x)\}_{j\in I} ~:=~\left\{ \sqrt{{\alpha_j\over k_j(\xi ^*)}}
P_j(x,\xi ^*)\right\}_{j\in I} ~\in ~P_I
\leqno (2.11)
$$
and using the definition  of ${\cal M}(\xi^*)$ 
and $\sum_{j\in {\cal M}(\xi^*)} \alpha_j = 1$ we obtain
$$
\sum_{j\in I} m_j^2(P_j^*)~=~\sum_{j\in {\cal M}(\xi ^*)}
m_j^2(P_j^*)~=~{1\over k_n(\xi^*)} ~=~\max \{ k_j^{-1}(\xi ^*)~|~
j\in I\}~.
$$
Therefore we have equality in Theorem 2.1 for $\{P_j^*\}_{j\in I}
\in P_I$ and $\xi^*\in \Xi$ and the assertion of the theorem follows.
\hfill \bull

\bigskip

{\bf Remark 2.3.} In the statistical theory the support points
and the weights of the optimal probability measure (minimizing 
(${\cal D}_I$)) give the relative frequencies and locations
of the observations in a polynomial regression. For the
special index set $I=\{1,\ldots ,n\}$ this measure has
been determined explicitly
 in Dette (1994).

\bigskip

{\bf Remark 2.4.} The polynomial $P_j^*$ in the set $\{
P_j^*\}_{j\in I}$ vanishes, whenever $j\notin {\cal M}(\xi^*)$
(which follows from $p_{2j}^*={1\over 2}$ and (2.9)), however, there
might be situations where $\alpha_j=0$ also for some $j\in {\cal M}(\xi ^*)$.
Observing the arguments at the end of the proof of the preceeding theorem
the solution of (${\cal P}_I$)  is obtained from (2.11) where the monic
polynomials (orthogonal with respect to the measure $d\xi^* (x)$) are
given by (2.4) and the quantities $k_j(\xi^*)$ are
obtained from (2.5). 
This provides a complete solution of the extremal problem (${\cal P}_I$).
In the following we will discuss some special cases of the set $I$
for which this solution  becomes more transparent.

\bigskip
\bigskip

{\bf 3. Chebyshev polynomials of the first kind.} If $I=\{n\}$, the solution
of $({\cal P}_I)$ is given by the Chebyshev polynomial of the first kind (on
the interval $[-b,b]$) $T_n({x\over b})$ (see Rivlin (1990) or  Natanson
(1955)). In this Section we will discuss two other sets for which the 
extremal polynomials have a  relative simple structure,
 namely $I =\{1,\ldots ,n\}$ and $I=\{n-1,n\}$. It turns out that the
answer of the question if the Chebyshev polynomial of the first kind is
also a solution of $({\cal P}_I)$ for these sets will depend heavily on the length
of the interval $[-b, b]$. We will start with the discussion of the problem
$({\cal P}_I)$ for the set $I=\{1,2,\ldots,n\}$. In the following $U_n(x)$ 
denotes the Chebyshev polynomial of the second kind (on the interval $[-1,1]$).

\bigskip

{\bf Theorem 3.1.}~~{\it Let $I=\{1,\ldots ,n\}$ and
$$
k~=~\min\Bigl\{ j\in \{1,\ldots ,n\}~|~U_{2n-2i+1}({b\over 2}) 
> 0 ~~\hbox{for } i=j,\ldots ,n \Bigr\}
\leqno (3.1)
$$
then the solution of the extremal problem $({\cal P}_I)$ 
is given by the polynomials
$\{P_l^*(x)\}_{l=1}^n$ where $P_l^*(x)=0$ if $l\leq k-1$,
$$
P_l^*(x)~=~\beta_l \left[ T_k({x\over b})U_{l-k}({x\over 2})~-~
{U_{n-k+1}({b\over 2})\over U_{n-k}({b\over 2})}
T_{k-1}({x\over b})U_{l-1-k}({x\over 2}) \right]
\leqno (3.2)
$$
($l=k,\ldots ,n$) and 
$$
\beta_l~=~ \pm {\sqrt{bU_{2n-2l+1}
({b\over 2})}\over U_{n-k+1}({b\over 2})} 
~~~~(l=k,\ldots ,n).
$$
The maximum value of $({\cal P}_I)$ is given by
$$ {2^{2k-2} \over b^{2k-1}} {U_{n-k}({b\over 2})\over
U_{n-k+1}({b\over 2})}.
$$
}
{\bf Proof.}~~For $j=n,\ldots ,k$ define $\gj=1-b^{-2(n-j)}
\prod_{i=j+1}^{n-1} (q_{2i}^*p_{2i}^*)^{-1}$
(here we put $\gamma_n(\xi^*) =1$
and the $p_{2j}^*$
are defined by (2.6)), then it is straightforward to show that
$$
\gj ~=~ {U_{n-j+1}({b\over 2}) \over b U_{n-j}({b\over 2})}~=~
{U_{2n-2j+1}({b\over 2})\over 2b U_{n-j}^2({b\over 2})}+{1\over 2}
~~~~~~~~j=k,\ldots ,n
\leqno (3.3)
$$
and the definition of $k$ in (3.1) and  Theorem 2.2 yield for the
canonical moments of the solution  $\xi^*$ of
the dual problem (${\cal D}_I$) $p_{2j}^*=\gamma_j (\xi^*)$
($j=k,\ldots ,n$). If $k\ge 2$, then it follows that  
$$
\gamma_{k-1}(\xi^*)~=~1-b^{-2(n-k+1)} \prod_{i=k}^{n-1}
(q_{2i}^*p_{2i}^*)^{-1} ~\leq ~{1\over 2}~
$$
 and that $b\leq 2$ which implies (by Theorem 2.2)
 $p_{2k-2}^*={1\over 2}$ and $\gamma_{k-2}  (\xi^*)
\leq 1-{2b^{-2}} 
\leq {1\over 2}$. Therefore the canonical moments of the solution $\xi^*$ 
of the dual problem (${\cal D}_I$) 
in Theorem 2.2 are given by 
$$
({1\over 2},{1\over 2},\ldots ,{1\over 2},
p_{2k}^*,{1\over 2},p_{2k+2}^*,{1\over 2},\ldots , {1\over 2},
p_{2n-2}^*,{1\over 2},1)
$$
where $p_{2j}^*= \gj$ ($j=k,\ldots ,n$) and $\gj$ is defined
in (3.3). By Theorem 2.1 we have to find the
orthonormal polynomials with respect to the measure
$d\xi^*(x)$ whose monic form is given by (2.4) that is
$$
\eqalign{
P_l(x,\xi^*)~&=~ K 
   \pmatrix {\overbrace{-{b^2\over 2} ~-{b^2\over 4} \ldots -{b^2
			\over 4}}^{k-1}& -{b^2\over 2} p_{2k}^*
~~~~-b^2q_{2k}^*p_{2k+2}^*& \ldots& -b^2q_{2l-4}^*p_{2l-2}^*
                       \cr
		  \noalign{\smallskip}
            x\quad\quad ~x\ \ \ldots\ \quad x& 
			{\qquad x \qquad\qquad\qquad x}& \ldots& 
			{x\hskip.8truein x}\cr}
 \cr
&=~ K 
   \pmatrix {\overbrace{-{b^2\over 2} ~-{b^2\over 4} \ldots -{b^2
			\over 4}}^{k-1}& -{b^2\over 2} p_{2k}^*
~~~~~~-1& \ldots& \! \! \! \! -1
                       \cr
		  \noalign{\smallskip}
            x\quad\quad ~x\ \ \ldots\ \quad x& 
			{\qquad \qquad x \qquad\qquad x}& \ldots& 
			{x\hskip.4truein x}\cr}
			\cr
			}
			$$
			$$
			=~2^{-k+1}b^k
\left[ T_k({x\over b})U_{l-k}({x\over 2})~-~
{U_{n-k+1}({b\over 2})\over U_{n-k}({b\over 2})}
T_{k-1}({x\over b})U_{l-1-k}({x\over 2}) \right].
\leqno (3.4)
$$
Here we have used Sylvester's identity (see e.g. Studden (1980), 
formula (4.12)), (3.3) (for $j=k$) and the recursive definition of the
Chebyshev polynomials of the first and second kind. Note that the case $k=1$ 
has to be considered separately but gives the corresponding result in (3.4) for $k=1$.
The $L_2$-norm of this polynomial is given by (note that 
$p_{2j}^*={1\over 2}$, $j=1,\ldots ,k-1$ and
$p_{2j}^*=\gj$, $j=k,\ldots ,n$)
$$
k_l(\xi^*) ~=~b^{2l}({1\over 4})^{k-1} p_{2k}^* \prod_{j=k+1}^l
q_{2j-2}^*p_{2j}^* ~=~ {b^{2k-1}\over 2^{2k-2}}
{U_{n-k+1}({b\over 2})\over U_{n-k}({b\over 2})}
$$
while the quantities $\alpha_l$ in (2.9) are obtained
as
$$\alpha_l~=~{U_{n-l}({b\over 2})[U_{n-l+1}({b\over 2}) - 
U_{n-l-1}({b\over 2})]
\over U_{n-k}({b\over 2})U_{n-k+1}({b\over 2})}
~=~ {U_{2n-2l+1}({b\over 2})\over U_{n-k}({b\over 2})
U_{n-k+1}({b\over 2})}
$$
($l=k,\ldots ,n$). The assertion now follows from Theorem 2.1.
\hfill \bull

\bigskip

{\bf 3.2 Discussion.}~~Theorem 3.1 shows that the structure
of the solution of (${\cal P}_I$) changes completely with
the length of the interval $[-b,b]$. If $b\leq \sqrt{2}$,
then we obtain from (3.1) $k=n$ and
consequently the sum of the squared leading
coefficients of the
polynomials $P_1^*,\ldots ,P_n^*$ is maximized for the choice
$P_l^*(x)=0$ ($1\leq l \leq n-1$) and $P_n^*(x)=
T_n({x\over b})$ with maximum value $(2^{n-1}b^{-n})^2$. If
$b>\sqrt{2}$ the situation changes completely. In this case the index 
$1\leq k \leq n$ defined by (3.1) depends on $n$ and $b$.
 The solution of the
problem (${\cal P}_I$) is given by (3.2). Finally, if $b\geq 2$, 
it follows that $k=1$ and (3.2) simplifies to
$$
P_l^*(x)~=~{\sqrt{U_{2n-2l+1}({b\over 2})}\over
\sqrt{b} U_n({b\over 2})}\left[
U_l({x\over 2}) ~-~{U_{n+1}({b\over 2})\over U_{n-1}({b\over 2})}
 U_{l-2}({x\over 2})\right]
~~~~l=1,\ldots ,n~.
$$

\bigskip

{\bf Example 3.3}~~Let $n=3$, then we have to distinguish the 
following cases:

A) If $b\leq \sqrt{2}$, we have $k=3$, the optimal polynomials
are given by
$$
P_1^*(x)=P_2^*(x)=0,~~P_3^*(x)=\pm T_3({x\over b})
$$
and the maximum is $16b^{-6}$.

B) If $\sqrt{2}\leq b\leq \sqrt{3}$, then $k=2$, the optimal
polynomials are
$$
\eqalign{
&P_1^*(x)=0,~~P_2^*(x)=\pm {b\sqrt{b^2-2}\over b^2-1} T_2({x\over b}) \cr
& P_3^*(x)= \pm {b\over b^2-1}\Bigl[{x} T_2({x\over b})-
{b^2-1\over b} T_1({x\over b})\Bigr] 
~=~\pm {1\over 2(b^2-1)}\Bigl[ b^2T_3({x\over b}) -
(b^2-2)T_1({x\over b}) \Bigr]\cr
}$$
and the maximum value is $4b^{-2}(b^2-1)^{-1}$.

C) If $b\ge \sqrt{3}$, then $k=1$, the optimal polynomials are
$$
\eqalign{
& P_1^*(x) =\pm {\sqrt{b^4-4b^2+3}\over  b^3-2b}x,~~~P_2^*(x)
=\pm {1\over b \sqrt{b^2-2}} \Bigl[U_2({x\over 2}) -{b^4-3b^2+1\over
b^2-1}\Bigr] \cr
& P_3^*(x)= \pm {1\over b^3-2b}\Bigl[U_3({x\over 2}) -{b^4-3b^2+1\over b^2-1} 
U_1({x\over 2}) \Bigr]  \cr}
$$
and the maximum value is $(b^2-1)/[b^2(b^2-2)]$.
\hfill \bull

\medskip

In the remaining part of this section we will consider the
index set $I=\{n-1,n\}$. Thus the problem is
to maximize the sum of the squared coefficients
$$
m_{n-1}^2(P_{n-1})~+~m_n^2(P_n)
\leqno (3.4)
$$
over the set of all polynomials  
(of degree $n-1$ and $n$) satisfying 
$$
P_{n-1}^2(x)~+~P_n^2(x)~\leq ~1 ~~~~\hbox{for all ~~} x\in [-b,b].
\leqno (3.5)
$$
The solution of this problem can be obtained by a similar
reasoning as in Theorem 3.1  for $k=n$ and $k=n-1$
and we omit the details in the proof of the following
result.

\bigskip

{\bf Theorem 3.4.}~~{\it The polynomials $P_{n-1}^*(x)$ and
$P_n^*(x)$ maximizing $(3.4)$ subject to the restriction 
$(3.5)$ are given by
$$
\eqalign{
& (P_{n-1}^*(x),P_n^*(x))~=~\Bigl(0,\pm T_n({x\over b})\Bigr) ~~~~~~~~~~~~~~~
~~~\hbox{if }~b\leq \sqrt{2},\cr
& \cr
& (P_{n-1}^*(x),P_n^*(x))~=~\Bigl(\pm {b\sqrt{b^2-2}
\over b^2-1}T_{n-1}({x\over b}), ~\pm {1\over 2(b^2-1)}
[b^2T_n({x\over b})-(b^2-2)T_{n-2}({x\over b})] \Bigr)
 \cr
}
$$
if $b\geq \sqrt{2}$. 
The maximum values in $(3.4)$ are given by $2^{2n-2}b^{-2n}$,
if $b\leq \sqrt{2}$, and by $2^{2n-4}b^{-(2n-4)}(b^2-1)^{-1}$
if $b\geq \sqrt{2}$, respectively.}

\bigskip

{\bf Remark 3.5.}~~For index sets of the form $I_m=
\{n-m+1,\ldots ,n\}$ the corresponding results are obtained 
similar to Theorem 3.4. The values of $b$ where the
structure of the
solution is changing, are obtained successively from (3.1) as
$b=\sqrt{2}$, $b=\sqrt{3}$, $b=\sqrt{2+\sqrt{2}}$,
$b= (\sqrt{5+\sqrt{5}})/\sqrt{2} ,\ldots $ (see also 
Example 3.3).

\bigskip
\bigskip

{\bf 4. Chebyshev polynomials of the second kind.}~~In this
section we will briefly discuss some generalizations
of the extremal properties of the Chebyshev polynomials
of the second kind. Let $I$ denote a subset of $\{0,1,\ldots ,n\}$
and define
$$
\tilde P_I~:=~\Bigl\{ (P_j)_{j\in I}~|~P_j \in \P_j, j\in I, ~
\sup_{x\in [-b,b]} (b^2-x^2)\sum_{j\in I} P_j^2(x) \leq 1\Bigr\}
$$
as the set of all polynomials $(P_j)_{j\in I}$
such that a weighted sup-norm of the sum of squares is
less or equal $1$. We are interested in 
the problem
$$
\max \Bigl\{ \sum_{l\in I} m_l^2(P_l) ~|~(P_l)_{l\in I} 
\in \tilde P_I \Bigr\}~~.
\leqno ( \tilde {\cal P}_I)
$$
If $I=\{n\}$ we obtain the well known extremal proerty of the Chebyshev
polynomials of the second kind $U_n(x)$, if $b=1$,
(see e.g. Achieser (1956), p. 250) and more generally 
of $U_n({x\over b})/b$, if $b>0$.
For the sake of brevity we will only state the 
generalizations corresponding
to the index sets $I=\{0,\ldots ,n\}$ and $I=\{n-1,n\}$. All
proofs can be obtained by a similar reasoning as in the
previous sections and are therefore omitted.

\bigskip

{\bf Theorem 4.1.}~~{\it Let $I=\{0,\ldots ,n\}$ and
$$
k~=~\min\Bigl\{ j\in \{0,\ldots ,n+1\}~|~{U_{2n-2i+3}({b\over 2}) 
} > 0 ~~\hbox{for } i=j,\ldots ,n+1 \Bigr\}
\leqno (4.1)
$$
then the solution of the problem $(\tilde {\cal P}_I)$ is given by the 
polynomials
$\{P_l^*(x)\}_{l=0}^n$ where $P_l^*(x)=0$ if $l\leq k-2$ and
$$
P_l^*(x)~=~\tilde \beta_l \left[ U_{k-1}({x\over b})U_{l-k+1}
({x\over 2})~-~
{U_{n-k+2}({b\over 2})\over U_{n-k+1}({b\over 2})}
U_{k-2}({x\over b})U_{l-k}({x\over 2}) \right]
\leqno (4.2)
$$
$(l=k-1,\ldots ,n)$, where 
$$
\tilde \beta_l~=~ \pm {\sqrt{U_{2n-2l+1}
({b\over 2})}\over \sqrt{b} U_{n-k+2}({b\over 2})} 
~~~~(l=k-1,\ldots ,n).
$$
The maximum value of $(\tilde {\cal P}_I)$ is given by
$$ {2^{2k-2} \over b^{2k-1}} {U_{n-k+1}({b\over 2})\over
U_{n-k+2}({b\over 2})}.
$$
}

\bigskip

{\bf Remark 4.2.}~~If $b\leq \sqrt{2}$ then it follows from (4.1)
that $k=n+1$ and the solution of $(\tilde {\cal P}_I)$
is given by the polynomials $P_l^*(x)=0$, $l=0,\ldots ,n-1$,
and $P_n^*(x)={1\over b}U_n({x\over b})$. As
in Discussion 3.2 it follows that for $b\geq 2$ we
have $k=1$ and the optimal polynomials are ``essentially'' 
independent of
the interval $[-b,b]$ and proportional to
the Chebyshev polynomials of the second kind, that is
$$
P_l^*(x)~=~\pm {\sqrt{U_{2n-2l+1}({b\over 2} )}\over 
\sqrt{b}U_{n+1}({b\over 2})}
U_l({x\over 2})~~~~~l=0,1,\ldots ,n~.
$$
with maximum value $U_n({b\over 2})[bU_{n+1}({b\over 2})]^{-1}$.
In the interval $[\sqrt{2},2]$ we have $1\leq k\leq n+1$
(depending on $b$ and $n$) and the solution of
($\tilde {\cal P}_I$) is given by (4.1) and (4.2).

\bigskip

{\bf Theorem 4.3.}~~{\it Let $b\leq \sqrt{2}$, then the solution
of the problem 
$$
\hbox{maximize}~~~~~m_{n-1}^2(P_{n-1})~+~m_n^2(P_n)
\leqno (4.3)
$$
subject to the restriction
$$
\sup_{x\in [-b,b] }(b^2-x^2)[P_{n-1}^2(x)+P_n^2(x)]~\leq ~1
\leqno (4.4)
$$
is given by the polynomials $P_{n-1}^*(x)=0$, $P_n^*(x)={1\over b}
U_n({x\over b})$ with optimum value $2^{2n}b^{-2n+2}$. If 
$b\geq \sqrt{2}$ the maximum in $(4.3)$ subject to $(4.4)$
is attained for the polynomials
$$
(P_{n-1}^*(x),P_n^*(x))~=~\Bigl(\pm {\sqrt{b^2-2}
\over b^2-1}U_{n-1}({x\over b}), ~\pm {b\over 2(b^2-1)}
[U_n({x\over b})-{(b^2-2)\over b^2}U_{n-2}({x\over b})] \Bigr)
$$
with maximum value $(2/b)^{2(n-1)}(b^2-1)^{-1}$.}

\bigskip
\bigskip

{\bf Appendix.} (Proof of Theorem 2.1)~~The proof of Theorem 2.1 
follows from a standard result in the theory of 
optimal design in mathematical statistics (see Pukelsheim (1993)).
To be precise let $I=\{i_1,\ldots ,i_m\}$,
$i=m+\sum_{j=1}^mi_j$, $f_j(x)=(1,x,\ldots ,x^j)'$ (where
$'$ denotes transposition) and define for a probability
measure $\xi$ on the interval $[-b,b]$
$$
M_j(\xi)~=~\int_{-b}^bf_j(x)f_j(x)'d\xi (x)~\in ~
\R^{(j+1)\times (j+1)} ~~~(j\in I)
$$
which is called moment matrix in the theory of optimal
design. In the
following we will collect all matrices $M_{i_1}(\xi),
\ldots ,M_{i_m}(\xi)$ in one big matrix
$$
M (\xi) = \pmatrix{M_{i_1} (\xi)& & \cr
					 & \ddots& \cr
					 & & M_{i_m} (\xi)\cr} 
   \in {\Bbb R}^{i \times i}
$$
and define two matrices by 
$$
K = \pmatrix{e_{i_1}& & \cr
		           & \ddots& \cr
		           & &  e_{i_m}\cr}
   \in {\Bbb R}^{i \times m}~~~
N = \pmatrix{N_{i_1} & & \cr
					 & \ddots& \cr
					 & & N_{i_m} \cr} 
   \in {\Bbb R}^{i \times i}
$$
where $e_j=(0,\ldots ,0,1)'\in \R^{j+1}$ 
is the $(j+1)$th unit vector ($j\in I$), $N_{i_j}$
are nonegative  $(i_j+1)\times (i_j+1)$ matrices 
(i.e. $N_{i_j} \geq 0$) and all
other entries in these matrices are $0$. Defining
$\Phi_{-\infty}(A)=\lambda_{\min}(A)$ where
$A\in \R^{m\times m}$, $A\geq 0$ and $\lambda_{\min}(A)$
denotes the minimum eigenvalue of $A$ we obtain  for
the polar function of $\Phi_{-\infty}$ (see Pukelsheim (1993),
p.149)
$\Phi_{-\infty}^{\infty}(A)=trace(A)$. By the duality theorem on
page 172 in the same reference it now follows that
(note that $k_j(\xi) =[e_j'M_j^{-1}(\xi)e_j]^{-1}$)
$$
\eqalign{
\max_{\xi \in \Xi} & \min\{ k_j(\xi)~|~j\in I\}~=~
\max_{\xi \in \Xi} \Phi_{-\infty}((K'M^{-1}(\xi)K)^{-1}) \cr
= ~& \min \Bigl\{[\Phi_{-\infty}^{\infty}(K'NK)]^{-1} ~|~ 
N\in \R^{i\times i},~
N \ge 0, ~ trace(M(\xi)N) \leq 1 ~~\forall \xi \in \Xi\Bigr\} \cr
= ~& \min \Bigl\{ (\sum\limits_{j\in I} e_j'N_je_j)^{-1}~|~N_j \in
\R^{(j+1) \times (j+1)},~N_j \geq 0 ~~
\forall j\in I,  \cr
&~~~~~~~~~~~~~~~~~~~~~~~~~~~~~~~~~~~~~~~~~~~
\sum_{j\in I} trace(M_j(\xi)N_j) \leq 1 ~~\forall \xi \in \Xi\Bigr\} \cr
= ~& \min \Bigl\{ (\sum_{j\in I} (e_j'a_j)^2)^{-1}~|~a_j \in \R^{j+1}
~~\forall j\in I, ~~
\sum_{j\in I} (f_j(x)'a_j)^2 \leq 1~~\forall x \in [-b,b]\Bigr\}\cr
= ~& \min \Bigl\{ (\sum_{j\in I} m_j^2(P_j))^{-1}~|~(P_j)_{j\in I} \in
P_I\Bigr\}.\cr
}
\leqno (A1)
$$
In order to go from the third to the fourth line in (A1) we have used
that
$$
\sum_{j\in I} trace(M_j(\xi)N_j)~=~
\sum_{j\in I} \int_{-b}^b f_j(x)'N_jf_j(x)d\xi (x) ~\leq ~1 ~~~~
\forall ~\xi \in \Xi
$$
is equivalent to the inequality  
$$
\sum_{j\in I} f_j(x)'N_jf_j(x) ~\leq ~1~~~~~~\forall x~\in [-b,b]
\leqno (A2)
$$
and the fact that the minimum value does not change if the
matrices $N_j$ are replaced by  matrices of the form $a_ja_j'$
(see the following discussion).
This proves the first part of the Theorem. For the second part we discuss equality in (A1) that is equality
in the duality theorem in Pukelsheim (1993) (p. 171,172) and
obtain
$$
\leqalignno{
&\sum_{j\in I} trace(M_j(\xi^*)N_j) ~=~1 & (A3) \cr
& M_j(\xi^*)N_j~=~{e_je_j'N_j\over e_j'M_j^{-1}(\xi^*)e_j}
~~~~~~~~~j\in I& (A4) \cr
& \min_{j\in I} \Bigl\{(e_j'M_j^{-1}(\xi^*)e_j)^{-1}
\Bigr\} ~\sum_{j\in I}e_j'N_je_j~=~
\sum_{j\in I} {e_j'N_je_j \over e_j'M_j^{-1}(\xi^*)e_j
}~=~1 .&
(A5) \cr
}
$$ 
Observing  that $k_j^{-1}(\xi^*)=e_j'M^{-1}_j(\xi^*)e_j$ ($j=1,\ldots ,n$)
we obtain by straightforward calculation as a solution of 
(A3) and (A4) $N_j=\alpha_j
a_ja_j'$ where $a_j =\sqrt{k_j(\xi^*)}M_j^{-1}(\xi^*)e_j$ ($j\in I$),
$\alpha_j \geq 0$ (because $N_j\geq 0$) and $\sum_{j\in I} \alpha_j
=1$. Finally it follows from (A5) that $\alpha_j=0$ whenever $j\notin {\cal M}
(\xi^*)$. By Corollary 2.3 in Dette (1994b)  the polynomials $P_l^*(x,\xi^*)=a_l'
f_l(x)$ are orthonormal with respect to the
measure $d\xi^* (x)$ which yield for the monic orthogonal polynomials
$P_l(x,\xi^*)=\sqrt{k_l(\xi^*)}a_l'f_l(x)$ ($l=1,\ldots ,n$).
Consequently a solution of the right hand side of (A1)
is given by $\{\sqrt{\alpha_j/k_j(\xi^*)}P_j(x,\xi^*)\}_{j\in I}$
where $P_j(x,\xi^*)$ is the $j$th monic orthogonal polynomial
with respect to the measure $d\xi^* (x)$ and 
 the $\alpha_j$ have to satisfy
$$
\sum_{j\in I} \alpha_j k_j^{-1}(\xi^*) P_j^2(x,\xi^*)~=~
\sum_{j\in I} f_j(x)'N_jf_j(x)~\leq 1 ~
$$
for all $x\in [-b,b]$. This completes the proof of Theorem 2.1.
\hfill \bull

\bigskip
\bigskip

\noindent
{\bf References}

{\parindent=0pt

N.I. Achieser, (1956), {\it Theory of Approximation}, Dover, New York.

 Chihara, T.S. (1978). {\it An Introduction to Orthogonal Polynomials},
Gordon and Breach, New York.

H. Dette (1994), Optimal designs for identifying the degree of a polynomial 
regression, Ann. Statist., to appear.

H. Dette (1994a), Extremal properties for ultraspherical polynomials,
Journal of Approxi- mation theory, {\bf 76}, 246--273.

H. Dette (1994b), New identities for orthogonal polynomials on a compact
interval, J. Math. Anal. Appl., {\bf 179}, 547--573.

 T.S. Lau, W.J. Studden (1988), On an extremal problem of
Fej\'er, { Journal of Approximation Theory\/}, {\bf 53},
184--194.

I.P. Natanson (1955), {\it Konstruktive Funktionentheorie}, Akademie Verlag, Berlin.

O. Perron (1954).
{\it Die Lehre von den Kettenbr\"uchen (Band I, II).}
B.G. Teubner, Stuttgart.

F. Pukelsheim (1993), {\it Optimal Design of Experiments}, Wiley, New York.

T.J. Rivlin (1990).
{\it Chebyshev polynomials.}
Wiley: New York.

W.J. Studden (1980).
$D_s$-optimal designs for polynomial regression using continued
fractions.
{\it Ann.\ Statist}., {\bf 8}, 1132--1141.

 W.J. Studden (1981), { On a problem of Chebyshev}, Journal of
 Approximation Theory, {\bf 29}, 253--260.

 H.S.Wall (1948), {\it Analytic theory of continued fractions}, Van Nostrand, 
New York.
}
\bye